\documentclass[a4paper,12pt]{amsart}
\usepackage{amssymb}
\usepackage{ifthen}
\usepackage{graphicx}
\nonstopmode
\newtheorem{thm}{Theorem}
\newtheorem{cor}{Corollary}
\newtheorem{lem}{Lemma}

\newtheorem{claim}{Claim}

\newtheorem{conj}{Conjecture}
\newtheorem{prob}{Problem}

\newtheorem{rem}{Remark}

\theoremstyle{definition}
\newtheorem{defn}{Definition}[section]
\newtheorem{example}{Example}

\newenvironment{pf}[1][]{%
 \vskip 1mm
 \noindent
 \ifthenelse{\equal{#1}{}}%
  {{\slshape Proof. }}%
  {{\slshape #1.} }%
 }%
{\qed\medskip}

\newcounter{alphabet}
\newcounter{tmp}
\newenvironment{Thm}[1][]{\refstepcounter{alphabet}%
\bigskip%
\noindent%
{\bf Theorem \Alph{alphabet}}%
\ifthenelse{\equal{#1}{}}{}{ (#1)}%
{\bf .} \itshape}{\vskip 8pt}

\makeatletter
\newcommand{\Ref}[1]{\@ifundefined{r@#1}{}{\setcounter{tmp}{\ref{#1}}\Alph{tmp}}}
\makeatother

\newenvironment{Lem}[1][]{\refstepcounter{alphabet}%
\bigskip%
\noindent%
{\bf Lemma \Alph{alphabet}}%
{\bf .} \itshape}{\vskip 8pt}

\newcommand{\IR}{{\mathbb R}}

\newcommand{\IC}{{\mathbb C}}
\newcommand{\ID}{{\mathbb D}}

\def\be{\begin{equation}}
\def\ee{\end{equation}}

\newcommand{\bee}{\begin{enumerate}}
\newcommand{\eee}{\end{enumerate}}

\newcommand{\blem}{\begin{lem}}
\newcommand{\elem}{\end{lem}}
\newcommand{\bthm}{\begin{thm}}
\newcommand{\ethm}{\end{thm}}
\newcommand{\bcor}{\begin{cor}}
\newcommand{\ecor}{\end{cor}}
\newcommand{\beg}{\begin{example}}
\newcommand{\eeg}{\end{example}}
\newcommand{\begs}{\begin{examples}}
\newcommand{\eegs}{\end{examples}}
\newcommand{\bdefe}{\begin{defn}}
\newcommand{\edefe}{\end{defn}}
\newcommand{\bprob}{\begin{prob}}
\newcommand{\eprob}{\end{prob}}
\newcommand{\bques}{\begin{ques}}
\newcommand{\eques}{\end{ques}}
\newcommand{\bei}{\begin{itemize}}
\newcommand{\eei}{\end{itemize}}

\newcommand{\bde}{\begin{deter}}
\newcommand{\ede}{\end{deter}}
\newcommand{\bca}{\begin{case}}
\newcommand{\eca}{\end{case}}
\newcommand{\bcl}{\begin{claim}}
\newcommand{\ecl}{\end{claim}}
\newcommand{\bcon}{\begin{conj}}
\newcommand{\econ}{\end{conj}}
\newcommand{\bcons}{\begin{conjs}}
\newcommand{\econs}{\end{conjs}}
\newcommand{\bprop}{\begin{propo}}
\newcommand{\eprop}{\end{propo}}
\newcommand{\br}{\begin{rem}}
\newcommand{\er}{\end{rem}}
\newcommand{\brs}{\begin{rems}}
\newcommand{\ers}{\end{rems}}
\newcommand{\bo}{\begin{obser}}
\newcommand{\eo}{\end{obser}}
\newcommand{\bos}{\begin{obsers}}
\newcommand{\eos}{\end{obsers}}
\newcommand{\bpf}{\begin{pf}}
\newcommand{\epf}{\end{pf}}
\newcommand{\ba}{\begin{array}}
\newcommand{\ea}{\end{array}}
\newcommand{\beq}{\begin{eqnarray}}
\newcommand{\beqq}{\begin{eqnarray*}}
\newcommand{\eeq}{\end{eqnarray}}
\newcommand{\eeqq}{\end{eqnarray*}}

\newcounter{minutes}
\divide\time by 60
\newcounter{hours}
\multiply\time by 60 \addtocounter{minutes}{-\time}

\begin{document}
\title[Injectivity of sections of convex harmonic mappings and convolution theorems]
{Injectivity of sections of convex harmonic mappings and convolution theorems}

\thanks{
File:~\jobname .tex,
          printed: \number\day-\number\month-\number\year,
          \thehours.\ifnum\theminutes<10{0}\fi\theminutes}

\author{Liulan Li and Saminathan Ponnusamy
}
\address{Liulan Li, Department of Mathematics and Computational Science,
Hengyang Normal University, Hengyang,  Hunan 421008, People's
Republic of China} \email{lanlimail2012@sina.cn}
\address{Saminathan  Ponnusamy,
Indian Statistical Institute (ISI), Chennai Centre,
SETS (Society for Electronic Transactions and security),
MGR Knowledge City, CIT Campus, Taramani,
Chennai 600 113, India.
}
\email{samy@isichennai.res.in, samy@iitm.ac.in}
\subjclass[2000]{Primary: 30C45}
\keywords{Harmonic mapping, partial sum, univalent, convex, starlike and close-to-convex mappings, harmonic convolution and
direction convexity preserving map.
}
\maketitle

\begin{abstract}
In the article the authors consider the class ${\mathcal H}_0$ of sense-preserving harmonic
functions $f=h+\overline{g}$ defined in the unit disk $|z|<1$ and normalized so that
$h(0)=0=h'(0)-1$ and $g(0)=0=g'(0)$, where $h$ and $g$ are analytic in the unit disk.
In the first part of the article we present two classes $\mathcal{P}_H^0(\alpha)$ and $\mathcal{G}_H^0(\beta)$ of functions from ${\mathcal H}_0$
and show that if $f\in \mathcal{P}_H^0(\alpha)$ and $F\in\mathcal{G}_H^0(\beta)$,
then the harmonic convolution is a univalent and close-to-convex harmonic function in the unit disk provided
certain conditions for parameters $\alpha$ and $\beta$ are satisfied. In the second part we study the harmonic sections (partial sums)
$$ s_{n, n}(f)(z)=s_n(h)(z)+\overline{s_n(g)(z)},
$$
where $f=h+\overline{g}\in  {\mathcal H}_0$, $s_n(h)$ and $s_n(g)$ denote the $n$-th partial sums of $h$ and $g$, respectively.
We prove, among others, that if $f=h+\overline{g}\in{\mathcal H}_0$ is a univalent harmonic convex mapping,  then $s_{n, n}(f)$ is
univalent and close-to-convex in the disk $|z|< 1/4$ for $n\geq 2$, and $s_{n, n}(f)$ is also convex in the disk  $|z|< 1/4$ for $n\geq2$ and $n\neq 3$.
Moreover,  we show that the section $s_{3,3}(f)$ of $f\in {\mathcal C}_H^0$ is  not convex in the disk $|z|<1/4$ but is
shown to be convex in a smaller disk.

\end{abstract}

\maketitle \pagestyle{myheadings} \markboth{L. Li and S. Ponnusamy}{Sections of convex harmonic mappings}

\section{Introduction and Main Results}

 One of the interesting features about univalent harmonic mappings $f$ is that if $f$ is convex
(resp. starlike, convex in a direction $\alpha$) in the unit disk ${\mathbb D}=\{z \in \IC:\, |z|<1\}$, then
it is not in general that the function $g$ defined by
$g(z)=r^{-1}f(rz)$ is convex (resp. starlike, convex in a direction $\alpha$), for $r<1$.
The aim of this article is to discuss properties such as convolution results and sections of univalent harmonic mappings in the plane. Our
theorems are generalization of known results for univalent analytic mappings which we now recall.


The class ${\mathcal S}$ of all univalent mappings $h$ analytic in  ${\mathbb D}$ normalized by $ h(0)=h'(0)-1=0$
is the central object in the study of univalent function theory, see \cite{Du-uni,Pomm}.
In 1928, Szeg\"o \cite{Sz} proved that if  $h\in {\mathcal S}$ then all sections $s_n(h)(z):=\sum_{k=1}^n a_k z^k$ of $h$ are univalent
in the disk $|z|<1/4$ and the number $1/4$ cannot be replaced by a larger one.
There exists considerable amount of results in the literature concerning sections of mappings from $\mathcal S$
and some of its various geometric subclasses mentioned later in this section.
 We refer the reader to  \cite[\S8.2, pp. 243--246]{Du-uni} for a general survey and to the recent
papers \cite{OS2,Obpo2,OS1,ObpoWir1} which stimulated further interest on this topic. Moreover, the
theory of Hadamard convolution also plays a major role in dealing with such problems. See \cite{FourSilver91,Goodman2,Rush88,silver88}.
However, corresponding questions for the class of univalent harmonic mappings seem to be  difficult to handle
as can be seen from the recent investigations of the authors \cite{Do,DN,LS,LS2}.

Let ${\mathcal H}$ be the class of all complex-valued harmonic functions $f=h+\overline{g}$ defined on ${\mathbb D}$,
where $h$ and $g$ are analytic on ${\mathbb D}$ with the normalization $h(0)=0=h'(0)-1$ and $g(0)=0$. Set
$${\mathcal H}_0 =\{f=h+\overline{g} \in {\mathcal H}:\, g'(0)=0\}.
$$
According to the work of Lewy \cite{lewy-36}, a function  $f=h+\overline{g}\in {\mathcal H}$ is locally univalent
and sense-preserving on $\ID$ if and only if its Jacobian $J_f(z)$ is positive in $\ID$, where
$$J_f(z)=|f_z(z)|^2-|f_{\overline z}(z)|^2=|h'(z)|^2-|g'(z)|^2.
$$
In view of this result, we observe that $J_f(z)>0$ in $\ID$ if and only if $h'(z)\neq 0$ in $\ID$
and the (second complex) dilatation  $\omega (z)=g'(z)/h'(z)$ of $f=h+\overline{g}$ is analytic in
$\ID$ and has the property that $|\omega (z)|<1$ for $z\in \ID$.

Following the pioneering work of  Clunie and  Sheil-Small \cite{CS}, let
${\mathcal S}_{H}$ denote the subclass of ${\mathcal H}$ that are sense-preserving
and univalent in $\ID$, and further let
${\mathcal S}_{H}^0={\mathcal S}_{H}\cap {\mathcal H}_0.
$
The class ${\mathcal S}_{H}^0$  reduces to ${\mathcal S}$ when  $g(z)$ is identically zero.
Note that each $f=h+\overline{g} \in {\mathcal H}_0$ has the form
\be\label{li5-eq1}
h(z)=z+\sum _{n=2}^{\infty}a_nz^n~\mbox{ and }~g(z)=\sum _{n=2}^{\infty}b_nz^n.
\ee
For $p\geq1$ and $q\geq2$, we define the harmonic {\em sections/partial sums} $s_{p,q}(f)$
of $f=h+\overline{g}\in {\mathcal H}_0$ as follows:
$$ s_{p, q}(f)(z)=s_p(h)(z)+\overline{s_q(g)(z)}.
$$
Also, denote by $\omega_{p,q}(f)$ the dilatation of the harmonic sections $s_{p,q}(f)(z)$.

Recall that a domain $\Omega$ is said to be close-to-convex if the complement of $\Omega$ can be
written as a union of non-intersecting half-lines. A harmonic function $f\in {\mathcal H}$
is said to be convex (resp. close-to-convex, starlike) in $|z|<r$ if it is univalent and
the range $f(|z|<r)$ is convex (resp. close-to-convex, starlike with respect to the origin).
By ${\mathcal C}_H^0$ (resp.  ${\mathcal K}_H^0$, ${\mathcal S}_H^{0*}$), we
denote the subclasses of functions in ${\mathcal S}_H^0$ which are convex (resp. close-to-convex, starlike)
in $|z|<1$ just like ${\mathcal C}$,  ${\mathcal K}$ and ${\mathcal S}^{*}$
are the subclasses of functions in ${\mathcal S}$ mapping $\ID$ onto these respective domains.
The reader is referred to \cite{CS,Du,SaRa2013} for many interesting results on planar univalent
harmonic mappings.

Szeg\"o \cite{Sz} also proved that if  $h\in {\mathcal C}$ (${\mathcal S}^*$), then all sections $s_n(h)$
of $h$ are convex (starlike) in the disk $|z|<1/4$.  Miki \cite{Mi} showed that the same holds for
close-to-convex functions in ${\mathcal S}$. We refer to
\cite{BhPo2014,iliev,Ma,OS2,OS1,samy-hiroshi-swadesh, robert41,Rush88,silver88,Si} for many
interesting results  and expositions on this topic in the case of conformal mappings.
In the case of univalent harmonic mappings, almost nothing is known in the literature until recently,
where for a given $\alpha<1$, the authors in  \cite{LS,LS2} considered the class
$$\mathcal{P}_H^0(\alpha)=\{f=h+\overline{g}\in\mathcal{H}_0:~{\rm Re}\,(h'(z)-\alpha)>|g'(z)| ~\mbox{for}~ z\in\ID\}
$$
and discussed properties of harmonic sections of functions from the class  ${\mathcal P}_H^0:=\mathcal{P}_H^0(0)$
(see Theorems \Ref{ThmA} and \Ref{ThmB}).
We note that functions in $\mathcal{P}_H^0(\alpha)$ are univalent and close-to-convex in the unit disk $\ID$ whenever $0\leq\alpha<1$.
Moreover, $\mathcal{P}_H^0(\alpha)\subset \mathcal{P}_H^0$ for $0\leq\alpha<1$ and
$ \mathcal{P}_H^0\subset \mathcal{K}_H^0$ so that ${\mathcal P}_H^0 \subsetneq {\mathcal S}_H^0$.
Also for $\beta<1$, we define
$$\mathcal{G}_H^0(\beta)=\{f=h+\overline{g}\in\mathcal{H}_0:\,{\rm Re}\left(\frac{h(z)}{z}\right)-\beta>\left|\frac{g(z)}{z}\right|
~\mbox{ for}~ z\in\ID\}
$$
and observe that $\mathcal{G}_H^0(\beta)\subset \mathcal{G}_H^0(0):=\mathcal{G}_H^0$ for $0\leq\beta<1$.
The classes $\mathcal{P}_H^0(\alpha)$ and $\mathcal{G}_H^0(\beta)$ will be considered to state and prove a new
convolution result (see Theorem \ref{BPL5-th1}) on the lines of ideas of Ponnusamy \cite{samy95} for analytic
functions.

We define the harmonic convolution (or Hadamard product) as follows:
For $f=h+\overline{g} \in {\mathcal H}$ with the series expansions for $h$ and $g$ as
in \eqref{li5-eq1}, and $F = H + \overline{G}\in {\mathcal H}$, where
$$
H(z)=z+\sum_{n=2}^\infty A_nz^n~\mbox{ and }~ G(z)=\sum_{n=1}^\infty B_nz^n,
$$
we define
$$(f\ast F)(z)=z+\sum_{n=2}^\infty a_nA_nz^n+\sum_{n=1}^\infty \overline{b_nB_n}\overline{z^n}.
$$
Clearly, $f\ast F=F\ast f$.  Thus, for two subsets $\mathcal{P}, \mathcal{Q}\subset \mathcal{H}$, we define
$\mathcal{P}*\mathcal{Q}=\{f*g:\, f\in\mathcal{P}, g\in\mathcal{Q}\}.$

\begin{thm}\label{BPL5-th1}
Let $\alpha, \beta$ $\in[0,1)$ and $\gamma =1-2(1-\alpha)(1-\beta)$. Then $\mathcal{P}_H^0(\alpha)*\mathcal{G}_H^0(\beta)
\subset\mathcal{K}_H^0$, whenever $\gamma\geq 0$.
In particular, $\mathcal{P}_H^0*\mathcal{G}_H^0(1/2) \subset \mathcal{K}_H^0$ and
$\mathcal{P}_H^0(1/2)*\mathcal{G}_H^0 \subset \mathcal{K}_H^0$.
\end{thm}

The proof of Theorem \ref{BPL5-th1} will be given in Section \ref{sec-conv1}. We now present
an example which shows that there are harmonic functions in $ \mathcal{G}_H^0(\beta)$ that are not univalent in $\ID$.

\beg\label{BPL5-eg1}
Consider the harmonic function $f(z)=z+a(1-\beta)\overline{z}^2$,  where $0\leq\beta<1$ and $a\in\mathbb{C}$. By the
definition of $\mathcal{G}_H^0 (\beta)$ it is clear that $f\in\mathcal{G}_H^0(\beta)$ if and only if $|a|\leq1$.  A
direct calculation shows that $f$ is univalent in $\ID$ if and only if $|a|\leq1/2(1-\beta)$. Thus if $a$ is a complex
number such that $|a|\in(1/2(1-\beta),1]$ then $f\in\mathcal{G}_H^0(\beta )$, but is not necessarily univalent in $\ID$.
\eeg

\br
Dorff  \cite{Do} (see also \cite{DN}) considered  ${\mathcal S}_H^0$ mappings that are convex in one direction and
these results have been extended by the present authors in \cite{LiPo1,LiPo2}. According to Theorem \ref{BPL5-th1}
and Example \ref{BPL5-eg1}, it follows that the convolution of a  non-univalent harmonic function with certain class of
harmonic functions could still be close-to-convex in $\ID$.  Note that $f(z)=z+(1/2)\overline{z}^2$ belongs to
$\mathcal{P}_H^0$ but is not convex in $\ID$.

At this place it is worth remarking the well-known fact that the convolution of two convex
functions in ${\mathcal C}_H^0$  is not necessarily univalent in $\ID$ (see also \cite{Do}). To do this, we consider
the harmonic convex mapping $f_0=h_0+\overline{g_0}  \in {\mathcal C}_H^0$, where
\be\label{eq-f_01}
h_0(z)=\frac{2z-z^2}{2(1-z)^2}~\mbox{ and }~ g_0(z)=\frac{-z^2}{2(1-z)^2}.
\ee
The function $f_0$ maps $\ID$ harmonically onto the half-plane $\{w:\,{\rm Re\,}w >-1/2\}$ and can be
obtained as the vertical shear (i.e. shear in the direction $\pi/2$) of the function $l(z)=z/(1-z)$ with dilatation
$\omega (z)=-z$. That is, $h_0$ and $g_0$ are obtained as the solution of the linear system
$$ h_0(z)+ g_0(z)=l(z)~\mbox{ and }~  g'_0(z)/ h'_0(z)=-z
$$
with the conditions $ h_0(0)= g_0(0)=0$  (see Shearing theorem due to Clunie and Sheil-Small \cite[Theorem 5.3]{CS}).
The function $f_0$ plays the role of extremal for certain extremal problems for the class ${\mathcal C}_{H}^0$.
Now, we see that the convolution $f_0*f_1$ of the right-half plane mapping
$f_0$ and the $6$-gon mapping (see \cite{Du-survey}) defined by $f_1=h_1+\overline{g_1}$, where
$$h_1(z)= z+\sum_{n=2}^\infty \frac{z^{6n+1}}{6n+1} ~\mbox{ and }~g_1(z)= -\sum_{n=2}^\infty \frac{z^{6n-1}}{6n-1},
$$
is not even locally univalent in $\ID$. This is because the dilatation $\omega _{f_0*f_1}$ of $f_0*f_1$ has the property that
$$|\omega _{f_0*f_1}(z)|=\left | \frac{(g_0*g_1)'(z)}{(h_0*h_1)'(z)}\right | =\left | \frac{z^4(2+z^6)}{1+2z^6}\right | \nless 1
~\mbox{ for every $z\in\ID$}.
$$
\er

In order to state other results, we need to recall some standard notations and results on harmonic mappings.

A domain $D\subset\mathbb{C}$ is said to be convex
in the direction $\alpha$ $(\alpha \in \IR)$ if for every $a\in \IC$ the set
$D\cap \{a+te^{i\alpha}:\, t\in \IR\}$ is either connected
or empty.  A univalent harmonic function
$f$ defined on $|z|<r$ is said to be {\it convex in the direction $\alpha$}
if $f(|z|<r)$ is convex in the direction $\alpha$.
We denote by ${\mathcal C}_H(\alpha)$ the family of normalized univalent harmonic functions
which are convex in the direction $\alpha$ in $\ID$.
We may set ${\mathcal C}_H^0(\alpha):={\mathcal C}_H(\alpha)\cap {\mathcal H}_0$.

Obviously, every function that is convex in the direction $\alpha$
$(0\leq \alpha< \pi)$ is necessarily close-to-convex, but the converse is
not true. Clearly, a convex function is convex in every direction. The class of
functions convex in one direction has been studied by many mathematicians
(see, for example, \cite{Do,HengSch70}) as a subclass of
functions introduced by Robertson \cite{Rober36}. The case $\alpha =0$ (resp. $\alpha =\pi/2$)
is referred to us convex in real (resp. vertical) direction.

\begin{figure}
\begin{center}
\includegraphics[height=6.0cm, width=5.5cm, scale=1]{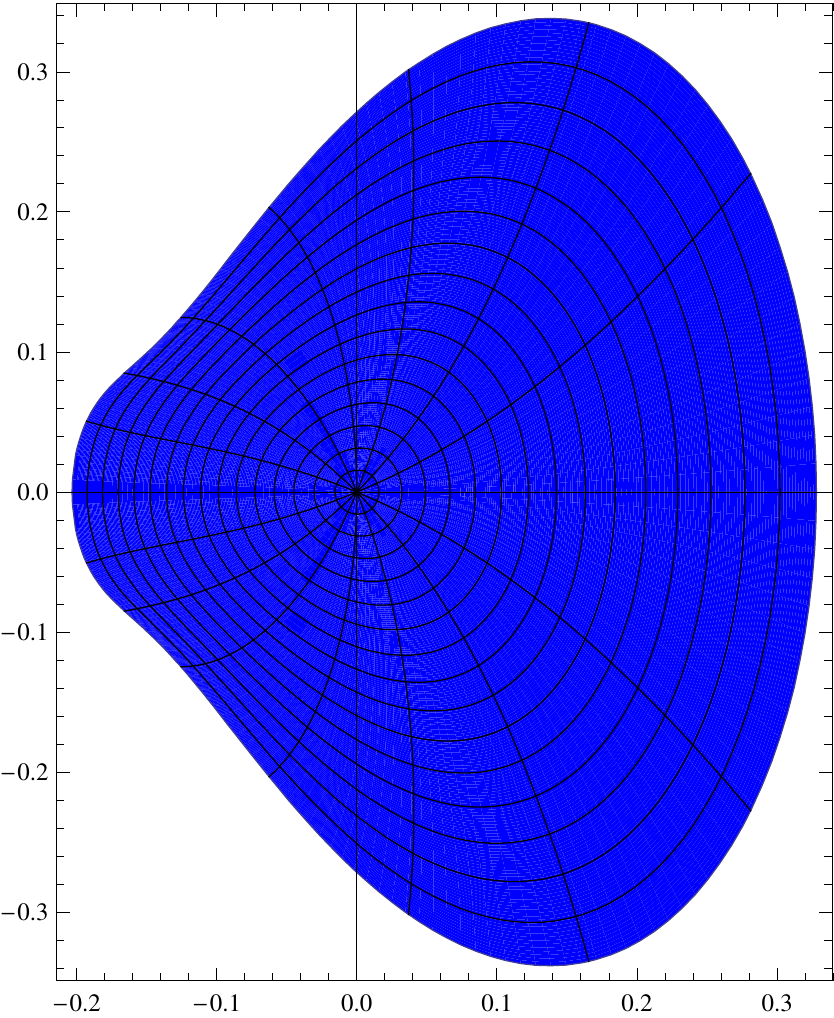}
\end{center}
\caption{Images of $\mathbb{D}_{1/4}$ under
$s_{3,3}(f_0)(z)=z + \frac{3 z^2}{2} + 2 z^3 - \overline{\big (\frac{z^2}{2} + z^3\big )}$\label{Fig1}}
\end{figure}

Concerning the classical result of Szeg\"o \cite{Sz} for the class
$\mathcal C$, it is natural to ask whether every section of $f\in {\mathcal C}_H^0 $ is convex
in some disk $|z|<r$. Thus, the first task is to derive properties of sections $s_{n,n}(f)$ of  $f\in {\mathcal C}_H^0 $.
Moreover, in our theorems we see that $s_{2,2}(f)$ and $s_{4,4}(f)$ are
(fully) convex in the disk $|z|<1/4$. It is surprising to see that
$s_{3,3}(f_0)$ is not convex in the disk $|z|<1/4$ (see Theorem \ref{33} and Figure \ref{Fig1}),
where $f_0$ is defined by \eqref{eq-f_01}.
This leads us to propose the following.

%

\begin{prob}\label{conj2}
Suppose that $f\in {\mathcal C}_{H}^0.$ Is each section $s_{n,n}(f)$ convex  in the disk $|z|<1/4$ for $n\geq2$ and $n\neq3$?
\end{prob}

In this article, we solve this problem and hence  our solution implies that for $n\geq2$ and $n\neq 3$, each section
$s_{n,n}(f)$ is convex in the direction of real axis in the disk $|z|<1/4$, in particular. On the other hand, Problem \ref{conj2}
remains open for the sections $s_{p,q}(f)$ of $f\in {\mathcal C}_{H}^0$ if $p\ne q$, $p\geq 1$ and $q\geq 2.$
Thus, as in the case of conformal mappings, it is natural to raise

\begin{prob}\label{prob2}
Suppose that $f\in {\mathcal S}_{H}^0$ (resp. ${\mathcal S}_{H}^{0*}, \, {\mathcal K}_{H}^{0},\, {\mathcal C}_{H}^{0},\,
   {\mathcal C}_{H}^{0}(\alpha))$.
Determine $\rho_{p,q}$ so that each section $s_{p,q}(f)$ belongs to the corresponding  class
in the disk $|z|<\rho_{p,q}$ for $p\geq 1$ and $q\geq 2$.
\end{prob}

Solution to Problem \ref{conj2} requires some ideas from the work of Ruscheweyh \cite{R} and, Ruscheweyh and Salinas \cite{RS}.


In Section  \ref{sec2}, we discuss the close-to-convexity of
$s_{n,n}(f)$. In Section \ref{sec3}, we prove that $s_{2,2}(f)$ of $f\in {\mathcal C}_{H}^0$ is
convex in the disk $|z|<1/4$ while $s_{3,3}(f_0)$ is not convex in
the disk $|z|<1/4$.  Finally,  in Section \ref{sec4}, we prove that (see Theorem \ref{main})
for $n\geq4$, each $s_{n,n}(f)$ is convex in the disk $|z|<1/4$.

We end this section with the following

\begin{conj}\label{conj1}
Suppose that $f\in {\mathcal C}_{H}^0.$ Then $s_{3,3}(f)$ is convex in the direction of real axis as well
as convex in the direction of imaginary axis, in the disk $|z|<1/4$.
\end{conj}

\section{Convolution theorem
}\label{sec-conv1}

We need  the following well-known result which follows easily from the Hergtlotz representation for
analytic functions with positive real part in the unit disk.

\begin{Lem}\label{BPL5-thA1}
If $p$ is analytic in $\ID$, $p(0)=1$, and ${\rm Re}\,p(z)>1/2$ in $\ID$ then for any function $F$, analytic in $\ID$, the
function $p*F$ takes values in the convex hull of the image of $\ID$ under
$F$.
\end{Lem}

We next recall another important result due to Clunie and Sheil-Small \cite{CS} which relates
the harmonic mapping $f=h+\overline{g}$ with the analytic functions $F_\lambda=h+\lambda g$.

\begin{Lem}\label{close-to-convexity}\cite{CS}
If a harmonic mapping $f=h+\overline{g}$ on $\ID$  satisfies $|g'(0)|<|h'(0)|$ and the
function $F_\lambda=h+\lambda g$ is close-to-convex for all
$|\lambda|=1$, then $f$ is close-to-convex and univalent in $\ID$.
\end{Lem}

%

\begin{pf}[Proof of Theorem \ref{BPL5-th1}]
Let $f_1\in\mathcal{P}_H^0(\alpha)$ have the canonical decomposition $f_1=h_1+\overline{g_1}$ with
\be\label{BPL5-eq1}
h_1(z)=z+\sum_{n=2}^{\infty}a_nz^n ~\mbox{and}~ g_1(z)=\sum_{n=2}^{\infty}b_nz^n.
\ee
Let $f_2\in\mathcal{G}_H^0(\beta)$ have the canonical decomposition $f_2=h_2+\overline{g_2}$ with
\be\label{BPL5-eq2}
h_2(z)=z+\sum_{n=2}^{\infty}A_nz^n ~\mbox{and}~ g_2(z)=\sum_{n=2}^{\infty}B_nz^n.
\ee
Now, we define $H=h_1*h_2+ \overline{g_1*g_2}$ and $H_{\epsilon}=(h_1*h_2)+\epsilon (g_1*g_2)$. Then
$H(0)=0=H_{\epsilon}(0)$ and $H_{\epsilon}'(0)=1$.
We need to show that $H\in\mathcal{K}_H^0$.
We remark that, as $(h_1*h_2)'(0)=1>(g_1*g_2)'(0)=0$, by Lemma \Ref{close-to-convexity},
it is enough to prove that, for all $\epsilon$ with $|\epsilon|=1$,
the function $H_{\epsilon}$ is close-to-convex in $\ID$.

By using the representations \eqref{BPL5-eq1} and \eqref{BPL5-eq2} we have
$$H_{\epsilon}'(z)= 1+\sum_{n=2}^{\infty}na_nA_nz^{n-1}+\epsilon\sum_{n=2}^{\infty}nb_nB_nz^{n-1}, ~\mbox{}~|\epsilon|=1.
$$
Now we claim that ${\rm Re}\,H_{\epsilon}'(z)>\gamma $, which will prove that $H_{\epsilon}$ is
in $\mathcal{P}_H^0(\gamma )$.

Since $f_1\in\mathcal{P}_H^0(\alpha)$,  the function $F_{\epsilon_1}$ defined by
$$F_{\epsilon_1}(z)= z+\frac{\sum_{n=2}^{\infty}a_nz^n+\epsilon_1(\sum_{n=2}^{\infty}b_nz^n)}{1-\alpha},~\mbox{}~z\in\ID,
$$
satisfies the condition ${\rm Re}\,F_{\epsilon_1}'(z)>0$,  for all $\epsilon_1$ with $|\epsilon_1|=1$.  A simple calculation shows that
the last inequality is equivalent to the inequality
\be\label{BPL5-eq3}
{\rm Re}\left(1+\frac{1}{2(1-\alpha)}\sum_{n=2}^{\infty}na_nz^{n-1}+\frac{\epsilon_1}
{2(1-\alpha)}\sum_{n=2}^{\infty}nb_nz^{n-1}\right)>\frac{1}{2},~\mbox{}~z\in\ID.
\ee
Similarly, as the function $f_2\in\mathcal{G}_H^0(\beta)$, for $|\epsilon_2|=1$ we have the inequality
$${\rm Re}\left(\frac{h_2(z)}{z}+\epsilon_2\frac{g_2(z)}{z}\right)>\beta,~\mbox{}~z\in\ID,
$$
which is equivalent to
\be\label{BPL5-eq4}
{\rm Re}\left(1+\frac{1}{2(1-\beta)}\sum_{n=2}^{\infty}A_nz^{n-1}+\frac{\epsilon_2}
{2(1-\beta)}\sum_{n=2}^{\infty}B_nz^{n-1}\right)>\frac{1}{2},~\mbox{}~z\in\ID.
\ee
Using Lemma \Ref{BPL5-thA1} and the inequalities \eqref{BPL5-eq3} and \eqref{BPL5-eq4} we get
$$
{\rm Re}\left(1+\frac{1}{4(1-\alpha)(1-\beta)}\sum_{n=2}^{\infty}na_nA_nz^{n-1}+\frac{\epsilon_1\epsilon_2}
{4(1-\alpha)(1-\beta)}\sum_{n=2}^{\infty}nb_nB_nz^{n-1}\right)>\frac{1}{2}.
$$
With $\gamma=1-2(1-\alpha)(1-\beta)$, the above inequality becomes
$${\rm Re}\left(1+\sum_{n=2}^{\infty}na_nA_nz^{n-1}+\epsilon_1\epsilon_2\sum_{n=2}^{\infty}nb_nB_nz^{n-1}\right)>\gamma ,~\mbox{}~z\in\ID,
$$
which shows that ${\rm Re}\, H_{\epsilon_1\epsilon_2}'(z)>\gamma $  for each $|\epsilon_1|=1$ and $|\epsilon_2|=1$.
In particular, for $\gamma\geq0$, $H_{\epsilon}(z)$ is close-to-convex for all $\epsilon$ with $|\epsilon|=1$. The proof is complete.
\end{pf}

%

\section{Close-to-convexity of sections $s_{n,n}(f)$ of convex functions $f$}\label{sec2}

%

By using Lemma \Ref{close-to-convexity} due to Clunie and Sheil-Small \cite{CS}, we obtain that

\begin{thm}\label{ctc}
Suppose that $f=h+\overline{g}\in {\mathcal H}_0$ is sense-preserving in $\ID$ and $F_\lambda=h+\lambda g$ is close-to-convex in $\ID$ for every
$|\lambda|=1$. Then $s_{n,n}(f)$ is close-to-convex and univalent in the disk $|z|<1/4$ for $n\geq2$.
\end{thm}
\bpf
Let $F_\lambda=h+\lambda g$ be close-to-convex. Then $f$ is locally univalent in $\ID$ and it follows that (see  Miki \cite{Mi}) $s_n(F_\lambda)$
is close-to-convex and univalent in the disk $|z|<1/4$ for all $n\geq2$. In other words, for each $n\geq2$,
the section  $4s_n(F_\lambda)(z/4)$ is close-to-convex and univalent in the unit disk $|z|<1$.
We observe that
$$4s_n(F_\lambda)(z/4)=4s_n(h)(z/4)+4\lambda s_n(g)(z/4),
$$
and so,
$$\left| \big (4s_n(h)(z/4)\big ) '(0)\right |=1>0=\left | \big (4s_n(g)(z/4)\big )'(0)\right |.
$$
By Lemma \Ref{close-to-convexity}, we find that
$$4s_n(h)(z/4)+\overline{4s_n(g)(z/4)}=4s_{n,n}(f)(z/4)
$$
is close-to-convex and univalent in the disk $|z|<1$ for all $n\geq2$.
The desired conclusion follows.
\epf


\br
We wish to emphasize that if  $f=h+\overline{g}\in{\mathcal S}_H^{0*}$, then it is not necessary that
the analytic functions $F_\lambda=h+\lambda g$ are univalent in $\ID$ for all $|\lambda|=1$.
For example, for $|\lambda |=1$, we consider
$$\varphi_\lambda(z)=\frac{z-\frac{1}{2}z^2+\frac{1}{6}z^3}{(1-z)^3}+\lambda \frac{\frac{1}{2}z^2+\frac{1}{6}z^3}{(1-z)^3}
=h(z)+\lambda g(z)
=z+\sum _{n=2}^{\infty}\varphi_{\lambda, n}z^n,
$$
where
$$\varphi_{\lambda, n} = \frac{1}{6} \left( 2n^2 (1+\lambda) + 3n(1-\lambda)+
(1+\lambda)\right) ~\mbox{ for all }~ n \geq 2.
$$
When $\lambda = -1$, $\varphi_\lambda (z)$ reduces to the analytic Koebe function $k(z)=z/(1-z)^2$, which is univalent
and starlike in $\mathbb{D}$. Moreover,  $\varphi_\lambda (z)$ is easily seen to be univalent only for $\lambda =-1$.
For $\varphi_\lambda$ to be univalent in $\ID$, it is necessary that $|\varphi_{\lambda,n}| \leq n$ for all $n \geq 2$.
For $|\lambda|=1$ ($\lambda \neq  -1$), we see that $|\varphi_{\lambda , n}| > n$ for large values of
$n$ and hence,  for these values of $\lambda$, $\varphi_\lambda (z)$ is not univalent in $\mathbb{D}$. Also, we observe that $K(z)=h+\overline{g}$ is
the harmonic Koebe mapping which is indeed starlike in $\ID$. This example shows that there is a limitation on the use of
Lemma \Ref{close-to-convexity}. However, analog of Theorem \ref{ctc} holds for the family ${\mathcal C}_H^0$ of univalent
harmonic convex mappings.
\er

\begin{thm}\label{locally univalent}
Let $f=h+\overline{g}\in{\mathcal C}_H^0$. Then every section $s_{n,n}(f)$ is close-to-convex
in the disk $|z|<1/4$ for $n\geq2$. In particular,
$s_{n,n}(f)$ is univalent and sense-preserving in $|z|<1/4$
for $n\geq2$. The number $1/4$ cannot be replaced by a grater one.
\end{thm}
\bpf
Let $f=h+\overline{g}\in{\mathcal C}_H^0$. Then the analytic functions $F_\lambda=h+\lambda g$
are close-to-convex in $\ID$ (see \cite[Theorem 5.7]{CS}) for all $|\lambda|=1$.
According to the last observation and
Theorem \ref{ctc}, we obtain that every section $s_{n,n}(f)$ is close-to-convex
in the disk $|z|<1/4$ for $n\geq2$.

Next we prove the sharpness part. Consider the function
$f_0=h_0+\overline{g_0}  \in {\mathcal C}_H^0$ defined by \eqref{eq-f_01}.
Then for $n=2$, we see that $s_2'(h_0)(z)=1+3z$ and $s_2'(g_0)(z)=-z$. Therefore, the dilatation
$\omega_{2,2}(f_0)$ of $f_0$ is given by
\be\label{li8-eq7}
\omega_{2,2}(f_0)(z)=\frac{s_2'(g_0)(z)}{s_2'(h_0)(z)}=\frac{-z}{1+3z}.
\ee
Since the M\"{o}bius transformation $w=M(z)=-z/(1+3z)$ maps the disk  $|z|<1/4$
onto the disk $|w-3/7|<4/7$, the relation \eqref{li8-eq7} implies that $|\omega_{2,2}(f_0)(z)|<1$ for $|z|<1/4$.
Moreover, at the boundary point $z=-1/4$, we have $\omega_{2,2}(f_0)(-1/4)=M(-1/4)=1$ which
shows that the radius $1/4$ cannot be replaced by a larger one.  The proof is complete.
\epf

\section{The sections $s_{2,2}(f)$ and $s_{3,3}(f_0)$}\label{sec3}

Let ${\mathcal A}_0$ denote the class of all functions $h(z)=\sum _{k=1}^{\infty}a_kz^k$ analytic
on the unit disk ${\mathbb D}$ and ${\mathcal A}=\{h\in  {\mathcal A}_0:\, h'(0)=1\}$.

A function $g\in{\mathcal A }_0$ is called \textit{Direction
Convexity Preserving} ($g\in {\rm DCP}$) if and only if $g\ast h\in
{\mathcal C}(\alpha)$ for all $h\in  {\mathcal C}(\alpha)$ and all $\alpha\in{\mathbb R}$.
Here ${\mathcal C}(\alpha)$ denotes the family of normalized univalent analytic functions in ${\mathbb D}$
which are convex in the direction $\alpha$.

The class ${\rm DCP}$ is somewhat special in the following sense: for
$g\in {\rm DCP}$, we do not necessarily have $g_{r}(z):=g(rz)\in {\rm DCP}$ for
$0<r<1$. We therefore define the ${\rm DCP}$ radius of an analytic
function $g$ to be $\max\{r:\, g_\rho\in {\rm DCP}\ \mbox{for}\
0<\rho<r\}$.

From \cite{R}, we observe that

\begin{Lem}\label{DCP}
$s_2(z)=z+z^2\in \rm{DCP}$ in the disk $|z|<1/4$.
\end{Lem}

We extend this lemma in Theorem \ref{DCP radius} for arbitrary section $s_n(z)$ of $z/(1-z).$
Let us now recall a convolution characterization for a function to be in the class $\rm{DCP}$.

\begin{Lem}\label{convexity preserving}\cite{RS}
Let $p\in {\mathcal A}_0$. Then $p\widetilde{\ast}f:=p\ast
h+\overline{p\ast g}\in {\mathcal C}_H^0$ for all
$f=h+\overline{g}\in {\mathcal C}_H^0$ if and only if $p\in {\rm DCP}$.
\end{Lem}

Before we proceed to state and prove our main results of this section, it is
appropriate to include the definition of (fully) convex mappings and
some known results on sections of functions from the class
${\mathcal P}_H^0$. For sense-preserving harmonic functions
$f=h+\overline{g} \in {\mathcal H}$, one has
$$\frac{\partial}{\partial \theta} \left (\arg\left
(\frac{\partial}{\partial \theta}f(re^{i \theta})\right )\right ) =
{\rm Re} \left (\frac{D^2f(z)}{Df(z)} \right )={\rm Re}\left
(\frac{z(h'(z) +zh''(z))+\overline{z(g'(z)+zg''(z))}}{zh'(z)-\overline{zg'(z)}}
\right ),
$$
where $z = re^{i \theta}$, $Df = z f_z - \overline{z}
f_{\overline{z}}$  and $D^2f = D(Df)$. Recall that if
$f=h+\overline{g} \in {\mathcal H}$ is sense-preserving,  $f(z)\neq
0 $ for all $z \in \mathbb{D} \backslash \{0\}$ and satisfies the
condition
$$ {\rm Re}\left
(\frac{z(h'(z) +zh''(z))+\overline{z(g'(z)+zg''(z))}}{zh'(z)-\overline{zg'(z)}}
\right )
 > 0  ~\mbox{ for all }~ z \in \mathbb{D} \backslash \{0\},
$$
then $f$ is univalent and \emph{fully convex} in $\ID$, i.e. the
image of every subdisk $|z|<r<1$ under $f$ is convex.

It is appropriate to recall two recent results of the authors.

\begin{Thm}{\rm (\cite[Theorems 4, 5 and 6]{LS})}\label{ThmA}
Let $f\in {\mathcal P}_H^0$. Suppose that $p$ and $q$ satisfy any
one of the following conditions:
\bee
\item[(a)] $p=1$ and $q\geq2$,
\item[(b)] $3\leq p<q$,
\item[(c)] $p=q\geq2$,
\item[(d)] $p>q\geq3$,
\item[(e)] $p=3$ and $q=2$.
\eee
Then $s_{p, q}(f)$ is univalent and close-to-convex in $|z|< 1/2$. Moreover, we have
\bee
\item[(f)] for $2<q$, $s_{2, q}(f)$ is univalent  and close-to-convex in $|z|<\frac{3-\sqrt{5}}{2}\approx 0.381966$.
\item [(g)] for $p\geq 4$, $s_{p, 2}(f)$ is univalent and close-to-convex in $|z|<0.433797$.
\eee
\end{Thm}

\begin{Thm} {\rm (\cite[Theorems 2, 3 and 4]{LS2})}\label{ThmB}
Let $f=h+\overline{g}\in {\mathcal P}_H^0$, and suppose that $p$ and $q$ satisfy one of the following
conditions:
\bee
\item[(a)] $p=1$ and  $q\geq2$
\item[(b)] $3\leq p<q$,
\item[(c)] $p=q\geq2$,
\item[(d)] $p>q\geq3$.
\eee
Then $s_{p, q}(f)$ is convex in $|z|<1/4$.
\bee
\item[(f)]
If $p=2<q$, then $s_{2, q}(f)$ is convex in $|z|<0.210222$.
\item[(g)] If $q=2<p$, then $s_{p, 2}(f)$ is convex in $|z|<0.234906$.
\eee
\end{Thm}

Now we explore the disk of convexity of $s_{n,n}(f)(z)$ when $f \in {\mathcal C}_H^0$. For $n=2$, we obtain

\begin{thm}\label{22}
Let $f=h+\overline{g}\in {\mathcal C}_H^0$, where $h(z)=z+\sum^\infty_{n=2}a_nz^n$ and
$g(z)=\sum^\infty_{n=2}b_nz^n$. Then the section $s_{2,2}(f)=z+a_2z^2+\overline{b_2z^2}$ is convex in
the disk $|z|<1/4$. The number $1/4$ cannot be replaced by a greater one.
\end{thm}
\bpf Set  $s_2(z)=z+z^2$. Then, by Lemmas \Ref{DCP} and \Ref{convexity preserving}, we conclude
that $r^{-1}s_2(rz)\widetilde{\ast}f(z)$ is convex in ${\mathbb D}$
for $0<r\leq\frac{1}{4}$. Since
$$r^{-1}s_2(rz)\widetilde{\ast}f(z)=z+ra_2z^2+\overline{rb_2z^2}=r^{-1}s_{2,2}(f)(rz),
$$
it follows that $r^{-1}s_{2,2}(f)(rz)$ is convex in ${\mathbb D}$
for $0<r\leq\frac{1}{4}$. This means that the section $s_{2,2}(f)$ is (fully) convex in
the disk $|z|<1/4$.

In order to prove the sharpness part, we consider the section $s_{2,2}(f_0)$ of $f_0=h_0+\overline{g_0}\in {\mathcal C}_H^0$,
where $h_0$ and $g_0$ are given by \eqref{eq-f_01}.
Note that
$$s_{2,2}(f_0)(z)=s_2(h_0)(z)+\overline{s_2(g_0)(z)}=z+(3/2)z^2-(1/2)\overline{z^2}.
$$
A computation gives
$${\rm Re}\left (\frac{z\left(zs_2'(h_0)(z)\right)'+\overline{z\left(zs_2'(g_0)(z)\right)'}}{zs_2'(h_0)(z)-\overline{zs_2'(g_0)(z)}}
\right ) ={\rm Re}\left (\frac{z+6z^2-2\overline{z}^2}{z+3z^2+\overline{z}^2}\right )={\rm
Re}\left (\frac{1+w(z)}{1-w(z)}\right ),
$$
where
$$w(z)=\frac{3z^2-3\overline{z}^2 }{2z+9z^2-\overline{z}^2} ~\mbox { and }~
\lim_{z\rightarrow0}\frac{1+w(z)}{1-w(z)} =1.
$$
Thus, for the convexity of $s_{2,2}(z)$ in the disk $|z|<1/4$, it
suffices to prove that $|w(z)|<1$ for $0<|z|<1/4$, which is equivalent
to
$$G(z)=\left|3z^2-3\overline{z}^2\right|^2-\left|2z+9z^2-\overline{z}^2\right|^2<0~\mbox { for }~0<\ |z|<1/4.
$$
Let $z=re^{i\theta}$. Then a computation yields
\beqq
G(re^{i\theta})&=&
36r^4\sin^2{2\theta}-\left[(2r\cos\theta+8r^2\cos{2\theta} )^2 +(2r\sin\theta+10r^2\sin{2\theta})^2\right ]\\
&=&36r^4\sin^2{2\theta}-\left(4r^2+64r^4+36r^4\sin^2{2\theta}+32r^3\cos\theta\cos{2\theta} +40r^3\sin\theta\sin{2\theta}\right)\\
&=& -\left [4r^2+64r^4 +32r^3\cos\theta (1-2\sin^2 \theta) + 64r^3\sin^2 \theta\cos \theta + 16r^3\sin^2\theta\cos{\theta}\right]\\
&=&-4r^2\left[1+16r^2+4r\cos\theta(2+\sin^2\theta)\right]\\
&=&-4r^2\left[1+16r^2+4r\cos\theta(3-\cos^2\theta)\right].
\eeqq
We observe that the function $B(x)=x(3-x^2)$ is increasing on $[-1,1]$
and therefore, from the last relation, we see that
$$G(re^{i\theta}) \leq -4r^2\left[1+16r^2+4rB(-1) \right]=-4r^2\left[1+16r^2-8r \right]
=-4r^2(4r-1)^2
$$
for $r<1/4$ and $-\pi <\theta\leq \pi$ with equality for $\theta =\pi$.
Thus, $G(z)<0$ for $0<|z|<1/4$ and hence, $|w(z)|<1$ for $|z|<1/4$. Finally,  $s_{2,2}(f_0)$ of $f_0=h_0+\overline{g_0}\in {\mathcal C}_H^0$
is (fully) convex for $|z|<1/4$ but not in a larger disk. The proof is complete.
\epf

For $n=3$, we will show that $s_{3,3}(f_0)(z)$ is not convex in
$|z|<1/4$.

\begin{thm}\label{33}
The harmonic section
$$s_{3,3}(f_0)(z)=s_3(h_0)(z)+\overline{s_3(g_0)(z)}=z+\frac{3}{2}z^2+2z^3-\frac{1}{2}\overline{z}^2-\overline{z}^3
$$
is not convex in the disk $|z|<1/4$. Here $f_0=h_0+\overline{g_0}\in {\mathcal C}_H^0$, where $h_0$ and $g_0$
are given by \eqref{eq-f_01}.
\end{thm}
\bpf By Theorem \ref{locally univalent}, $s_{3,3}(f_0)(z)$ is
locally one-to-one and sense-preserving in $|z|<1/4$. Now, by a computation, we have
\be\label{eq-main1a}
F(z)={\rm Re}\left (\frac{z\left(zs_3'(h_0)(z)\right)'+\overline{z\left(zs_3'(g_0)(z)\right)'}}
{zs_3'(h_0)(z)-\overline{zs_3'(g_0)(z)}} \right ) = {\rm Re}\left
(\frac{z+6z^2+18z^3-2\overline{z}^2-9\overline{z}^3}{z+3z^2+6z^3+\overline{z}^2+3\overline{z}^3}\right
).
\ee
Let $z_0=\frac{1}{4}e^{\frac{2i\pi}{3}}$. Then, it follows that
\beqq
F(z_0)&=&{\rm Re}\left
(\frac{\frac{1}{4}e^{\frac{2i\pi}{3}}+\frac{6}{16}e^{\frac{4i\pi}{3}}-\frac{2}{16}e^{\frac{2i\pi}{3}}+\frac{9}{64}}
{\frac{1}{4}e^{\frac{2i\pi}{3}}+\frac{3}{16}e^{\frac{4i\pi}{3}}+\frac{1}{16}e^{\frac{2i\pi}{3}}+\frac{9}{64}}\right )\\
&=&{\rm Re}\left
(\frac{\frac{1}{8}e^{\frac{2i\pi}{3}}+\frac{3}{8}e^{\frac{-2i\pi}{3}}+\frac{9}{64}}
{\frac{5}{16}e^{\frac{2i\pi}{3}}+\frac{3}{16}e^{\frac{-2i\pi}{3}}+\frac{9}{64}}\right )\\
&=&{\rm Re}\left (\frac{-\frac{7}{64}-\frac{\sqrt{3}}{8}i}{-\frac{7}{64}+\frac{\sqrt{3}}{16}i}\right )
=-\frac{47}{97}<0.
\eeqq
This means that $s_{3,3}(f_0)(z)$ is not convex in the disk $|z|<1/4$.
 \epf

\br
For the function $f_0=h_0+\overline{g_0}\in{\mathcal C}_H^0$ defined  by \eqref{eq-f_01}, it can be easily seen
that the function $F(z)$ defined by \eqref{eq-main1a} satisfies the positivity condition
$F(z)>0$ for $|z|<0.201254$ and thus, the disk of convexity of $s_{3,3}(f_0)$ is $|z|<r$, where $r$ is close to the
value $0.201254$. Since the computation is lengthy, we do not wish to address it for the moment. However,
in Theorem \ref{main1a}, we actually show that the section $s_{3,3}(f)(z)$ of every $f=h+\overline{g}\in{\mathcal C}_H^0$
is indeed convex in the disk $|z|<0.201254$.
\er

\begin{figure}
\begin{center}
\includegraphics[height=6.0cm, width=5.5cm, scale=1]{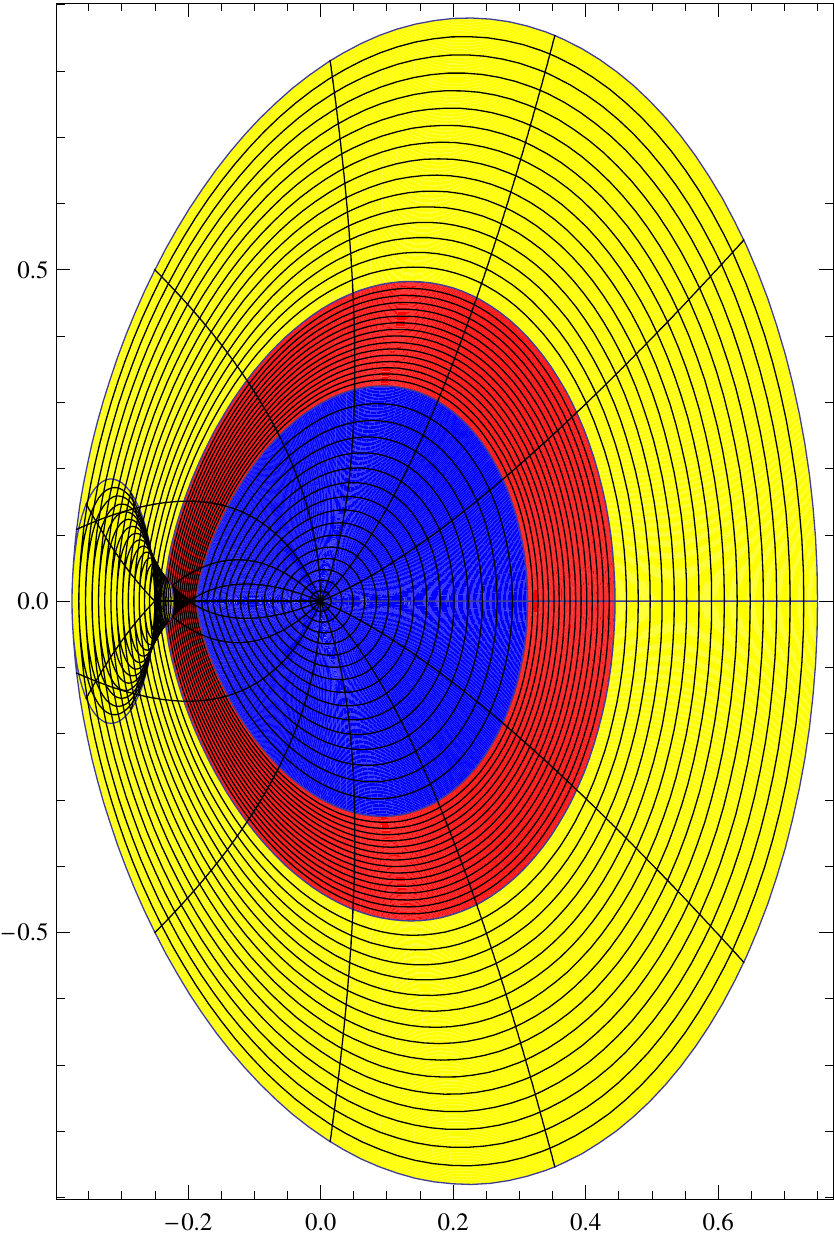}
\hspace{1cm}
\includegraphics[height=6.0cm, width=5.5cm, scale=1]{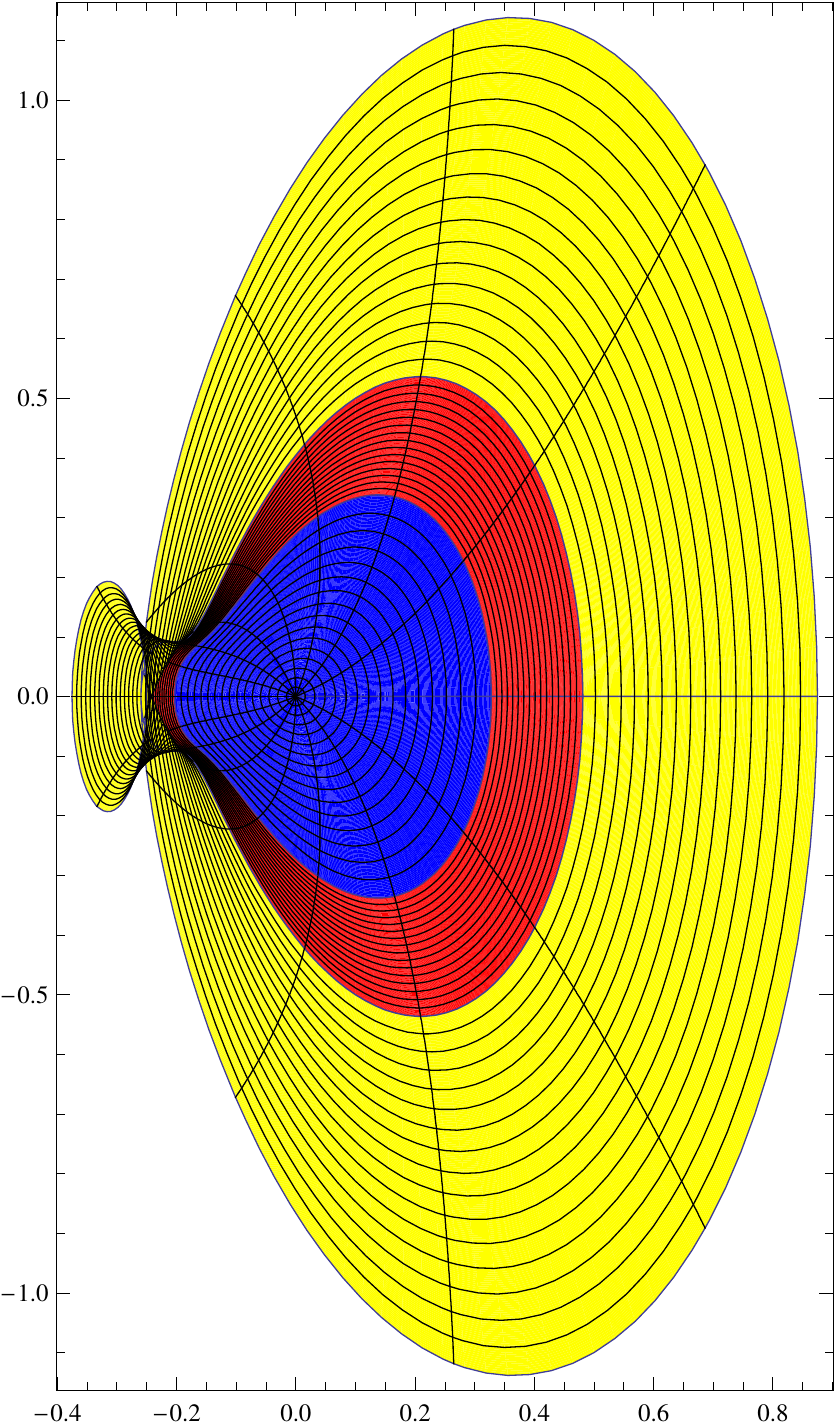}
\end{center}
\caption{Images of $|z| < 1/4$, $1/4 < |z| < 1/3$ and $1/3 < |z| <
1/2$  under $s_{2,2}(f_0)(z)$ and $s_{3,3}(f_0)(z)$
\label{PS4_Fig1}}
\end{figure}

In Figure \ref{PS4_Fig1}, images of $|z| < 1/4$, $1/4 < |z| < 1/3$
and $1/3 < |z| < 1/2$ under $s_{2,2}(f_0)(z)$ and $s_{3,3}(f_0)(z)$
are drawn in blue, red and yellow colors, respectively. These
pictures were drawn using  Mathematica as plots of the images of equally spaced
radial segments and concentric circles of the corresponding disk and of the two annuli.

\section{Disk of convexity of $s_{n,n}(f)$}\label{sec4}

We need the following result for the proof of two remaining theorems.

\begin{Lem}\label{DCP criterion2}\cite[Theorem 2]{RS}
Let $g$ be analytic in $\mathbb{D}$. Then $g\in {\rm DCP}$ if and only if
for each $t\in\mathbb{R}$, $g+itzg'$ is convex in the direction of
imaginary axis.
\end{Lem}

For the proof of Theorem \ref{DCP radius},   we use
a result of Royster and Ziegler \cite{RZ} concerning analytic mappings convex in one direction.



\begin{Lem}\label{lem3.2}\cite[Theorem 1]{RZ} Let $\phi(z)$ be a non-constant function analytic in $\mathbb{D}$. The
function $\phi(z)$ maps  univalently $\mathbb{D}$ onto a domain
convex in the direction of imaginary axis if and only if there are
numbers $\mu$ and $\nu$, $0\leq \mu <2\pi$ and $0\leq \nu\leq \pi$,
such that
\be\label{RZ-eq1}
{\rm Re\,}\{F_{\mu,\, \nu}(z)\phi'(z)\}\geq 0
\ee
for all $z\in \mathbb{D}$,  where $F_{\mu,\, \nu}(z)=-ie^{i\mu}(1-2ze^{-i\mu}\cos \nu+z^2e^{-2i\mu})$.
\end{Lem}

By using Lemmas \Ref{DCP criterion2} and \Ref{lem3.2}, we now prove the following theorem for $n\geq 4$
and in view of the technical details we present the proof of the case $n=3$ separately in Theorem \ref{main1a}

\begin{thm}\label{DCP radius}
$s_n(z):=\sum^n_{k=1}z^k=\frac{z-z^{n+1}}{1-z}\in {\rm DCP}$ in the disk $|z|<1/4$ for $n\geq4$.
\end{thm}
\bpf Let $\phi(z)=s_n(z)+itzs'_n(z)$, where $t\in\mathbb{R}$. A
computation yields that
\beqq
\phi'(z)&=&\frac{1-(n+1)z^n+nz^{n+1}}{(1-z)^2}+it\frac{1-(n+1)^2z^n+n(n+2)z^{n+1}}{(1-z)^2}\\
&&-it\frac{2nz^{n+1}}{(1-z)^2}+it\frac{2\sum^n_{k=1}z^k}{(1-z)^2}\\
&=&\frac{1-(n+1)z^n+nz^{n+1}}{(1-z)^2}+it\frac{1-(n+1)^2z^n+n^2z^{n+1}+2\sum^n_{k=1}z^k}{(1-z)^2}.
\eeqq
We now divide our proof into the following three cases.
\medskip

\noindent{\bf Case 1: $t>\frac{2}{19}$.}

Let $\mu=\nu=0$. Then $F_{0,\, 0}(z)=-i(1-z)^2$. It follows that
$$F_{0,\, 0}(z)\phi'(z)=t\left[1-(n+1)^2z^n+n^2z^{n+1}+2\sum^n_{k=1}z^k\right]-i\left[1-(n+1)z^n+nz^{n+1}\right]
$$
and
$${\rm Re\,}\{F_{0,\, 0}(z)\phi'(z)\}\geq
t-t\left[(n+1)^2|z|^n+n^2|z|^{n+1}+2\sum^n_{k=1}|z|^k\right]-(n+1)|z|^n-n|z|^{n+1}.
$$
It suffices to prove that the right side of the above inequality is
larger than $0$ for $|z|=\frac{1}{4}$ and for all $n\geq4$, since it is
harmonic in $|z|<\frac{1}{4}$. For $|z|=\frac{1}{4}$, the above
estimate takes the following form
\beqq
{\rm Re\,}\{ F_{0,\, 0}(z)\phi'(z)\}&\geq&
t-t\left[\frac{5n^2+8n+4}{4^{n+1}}+\frac{2-\frac{2}{4^n}}{3}\right]-\frac{5n+4}{4^{n+1}}\\
&=&\frac{t}{3}-\frac{5tn^2+(8t+5)n+\frac{4t}{3}+4}{4^{n+1}}:=A(n).
\eeqq
We see that $A(n)$ is monotonically increasing with respect to $n$ for $n\geq4$.
It follows that
$$A(n)\geq A(4)=\frac{57t}{4^4}-\frac{6}{4^4}=\frac{3}{4^4}(19t-2)>0
$$
for $t>\frac{2}{19}$, which implies that ${\rm Re\,}\{F_{0,\, 0}(z)\phi'(z)\}>0$ for
$n\geq4$ and $|z|=\frac{1}{4}$. Lemma \Ref{lem3.2} implies that
$\phi(z)$ is convex in the direction of imaginary axis in the disk
$|z|<1/4$ if $t>\frac{2}{19}$ and $n\geq4$.

\medskip

\noindent{\bf Case 2: $t<-\frac{2}{19}$.}

Let $\mu=\nu=\pi$. Then $F_{\pi,\, \pi}(z)=i(1-z)^2$. It follows that
$$F_{\pi,\, \pi}(z)\phi'(z)=-t\left[1-(n+1)^2z^n+n^2z^{n+1}+2\sum^n_{k=1}z^k\right]+i\left[1-(n+1)z^n+nz^{n+1}\right].
$$
By a similar reasoning as in Case $1$, we obtain that ${\rm Re\,}\{F_{\pi,\,\pi}(z)\phi'(z)\}>0$ for $n\geq4$ and
$|z|<\frac{1}{4}$. By Lemma \Ref{lem3.2}, we thus have shown that
$\phi(z)$ is convex in the direction of imaginary axis in the disk
$|z|<1/4$ if $t<-\frac{2}{19}$ and $n\geq4$.
\medskip

\noindent{\bf Case 3: $-\frac{2}{19}\leq t\leq \frac{2}{19}$.}

Let $\mu=\nu=\frac{\pi}{2}$. Then $F_{\pi/2,\, \pi/2}(z)=1-z^2=(1-z)(1+z).$ It follows that
\beqq
F_{\pi/2,\, \pi/2}(z)\phi'(z)&=&\frac{1+z}{1-z}+\frac{1+z}{1-z}\left(nz^{n+1}-(n+1)z^n\right)
+it\frac{1+z}{1-z}\\
&& +it\frac{1+z}{1-z}\left(-(n+1)^2z^n+n^2z^{n+1}+2\sum^n_{k=1}z^k\right),
\eeqq
and therefore,
\beqq
{\rm Re\,}\left (F_{\pi/2,\, \pi/2}(z)\phi'(z)\right )
&\geq&\frac{1-|z|}{1+|z|}-\frac{2|t|\cdot|z|}{(1-|z|)^2}-\frac{1+|z|}{1-|z|}\left(n|z|^{n+1}+(n+1)|z|^n\right)\\
&& -|t|\frac{1+|z|}{1-|z|}\left((n+1)^2|z|^n+n^2|z|^{n+1}+2\sum^n_{k=1}|z|^k\right).
\eeqq
For $|z|=\frac{1}{4}$, the above estimate takes the following form
\beqq
{\rm Re\,}\left (F_{\pi/2,\, \pi/2}(z) \phi'(z) \right )
&\geq&\frac{3}{5}-\frac{8|t|}{9}-\frac{5}{3}\cdot\frac{5n+4}{4^{n+1}}-
\frac{5|t|}{3}\left(\frac{5n^2+8n+4}{4^{n+1}}+\frac{2}{3}-\frac{2}{3}\cdot\frac{1}{4^n}\right)\\
&=&\frac{3}{5}-\frac{18|t|}{9}-\frac{5}{3}\cdot\frac{5|t|n^2+(8|t|+5)n+4+\frac{4|t|}{3}}{4^{n+1}}:=B(n).
\eeqq
We observe that $B(n)$ is monotonically increasing with respect to $n$ for $n\geq4$.
Hence,
$$B(n)\geq B(4)=\frac{1}{4^3}\left(\frac{359}{10}-\frac{5033|t|}{36}\right)>0
$$
for $-\frac{2}{19}\leq t\leq \frac{2}{19}$. Again, by Lemma \Ref{lem3.2},
we obtain that  $\phi(z)$ is convex in the direction of imaginary axis in the
disk $|z|<1/4$ if $-\frac{2}{19}\leq t\leq \frac{2}{19}$ and for all
$n\geq4$.
The  desired conclusion follows from Lemma \Ref{DCP
criterion2}.
 \epf

\begin{thm}\label{main}
Let $f=h+\overline{g}\in{\mathcal C}_H^0$. Then $s_{n,n}(f)$ is
convex  in the disk $|z|<1/4$ for $n\geq4$.
\end{thm}
\bpf By Theorem \ref{DCP radius} and Lemma \Ref{convexity
preserving}, we conclude that
$r^{-1}s_n(rz)\widetilde{\ast}f(z)$ is convex in ${\mathbb D}$ for
$0<r\leq\frac{1}{4}$ and $n\geq4$. Since
$$r^{-1}s_n(rz)\widetilde{\ast}f(z)=r^{-1}s_{n,n}(f)(rz),
$$
it follows that $r^{-1}s_{n,n}(f)(rz)$ is convex in ${\mathbb D}$
for $0<r\leq\frac{1}{4}$ and $n\geq4$. This means that the section
$s_{n,n}(f)$ is (fully) convex in the disk $|z|<1/4$ for $n\geq4$.
\epf




\begin{thm}\label{main1a}
Let $f=h+\overline{g}\in{\mathcal C}_H^0$. Then $s_{3,3}(f)$ is
convex  in the disk $|z|<0.201254$.
\end{thm}
\bpf
As in the proof of Theorem \ref{main}, it suffices to show that  $s_3(z):=z+z^2+z^3\in {\rm DCP}$ in the disk $|z|<0.201254$.

We only have to give the crucial steps and appropriate replacements in the proof of Theorem \ref{DCP radius} with $n=3$
and the rest of arguments follows from there. Thus, if $\phi(z)$ is as in the proof of Theorem \ref{DCP radius} with $n=3$, then
$\phi'(z)$ takes the form
$$\phi'(z)=\frac{1-4z^3+3z^4}{(1-z)^2}+it\frac{1-16z^3+9z^4+2\sum^3_{k=1}z^k}{(1-z)^2}.
$$
\noindent{\bf Case 1: $t>0.105712$.}

It follows from the proof of Theorem \ref{DCP radius} that
$${\rm Re\,}\{F_{0,\ 0}(z)\phi'(z)\}\geq t-t\left[16|z|^3+9|z|^4+2\sum^3_{k=1}|z|^k\right]-4|z|^3-3|z|^4,
$$
which for $|z|\leq 0.201254$ implies that
\beqq
{\rm Re\,}\{ F_{0,\,0}(z)\phi'(z)\}&\geq& t\left[1 - 16 (0.201254)^3 - 9
(0.201254)^4-2\sum^3_{k=1}(0.201254)^k\right]\\
&& -4(0.201254)^3-3(0.201254)^3>0
\eeqq
for $t>t_0 \approx 0.10571184$.
In particular, by Lemma \Ref{lem3.2}, we obtain that  $\phi(z)$ is convex in the direction of
imaginary axis in the disk $|z|<0.201254$ if $t>0.105712$.

\medskip

\noindent{\bf Case 2: $t<-0.105712$.}

With  $\mu=\nu=\pi$ so that  $F_{\pi,\, \pi}(z)=i(1-z)^2$, we have
$$F_{\pi,\, \pi}(z)\phi'(z)=-t\left[1-16z^3+9z^4+2\sum^3_{k=1}z^k\right]+i\left[1-4z^3+3z^4\right]
$$
and by a similar reasoning as in Case $1$, we obtain that
$${\rm Re\,}\{F_{\pi,\, \pi}(z)\phi'(z)\}>0 ~\mbox{ for $|z|<0.201254$}
$$
and thus, $\phi(z)$ is convex in the direction of imaginary axis in the disk $|z|<0.201254$ if
$t<-0.105712$.
\medskip

\noindent{\bf Case 3: $-0.105712\leq t\leq 0.105712$.}

This case corresponds to $\mu=\nu=\frac{\pi}{2}$ so that  $F_{\pi/2,\, \pi/2}(z)=1-z^2 $ and
\beqq
{\rm Re\,}\left (F_{\pi/2,\, \pi/2}(z)\phi'(z)\right )
&\geq&\frac{1-|z|}{1+|z|}-\frac{2|t|\cdot|z|}{(1-|z|)^2}-\frac{1+|z|}{1-|z|}\left(3|z|^4+4|z|^3\right)\\
&& -|t|\frac{1+|z|}{1-|z|}\left(16|z|^3+9|z|^4+2\sum^3_{k=1}|z|^k\right).
\eeqq
For $|z|=0.201254$, the above estimate shows that
$${\rm Re\,}\left (F_{\pi/2,\, \pi/2}(z)\phi'(z)\right )\geq 0.608489-1.60093|t|>0
$$
for $|t|<0.608489/1.60093 ~(> 0.105712)$.
Consequently, by Lemma \Ref{lem3.2}, we obtain that  $\phi(z)$ is convex in the direction of imaginary axis
in the disk $|z|<0.201254$ if $-0.105712\leq t\leq 0.105712$ and for $n=3$.

The cases $1$ to $3$ show that $s_3(z):=z+z^2+z^3\in {\rm DCP}$ in the disk $|z|<0.201254$.
%
%
\epf

\subsection*{Acknowledgements} The visit and the research of first author was supported by
``Abel visiting Scholar Program" of Commission for Developing Countries (IMU).
The research was also supported by NSF of China (No. 11201130),
Hunan Provincial Natural Science Foundation of China (No. 14JJ1012),
Scientific Research Fund of Hunan Provincial Education Department (No. 11B019) and
construct program of the key discipline in Hunan province.
The second author is currently on leave from the Indian Institute of Technology Madras.

\end{document}